\documentclass[twoside,12pt,a4paper]{amsart}

\usepackage{amsmath}
\numberwithin{equation}{section}

\begin{document}

\title {An Alternative to Moment Closure}

\author{Ingemar N{\aa}sell}
\address{Department of Mathematics \\
        The Royal Institute of Technology \\
        S-100 44 Stockholm, Sweden}
\email{ingemar@kth.se}
 

\maketitle

\begin{abstract}
Moment closure methods are widely used to analyze mathematical models. 
They are specifically geared toward derivation of approximations of moments of stochastic models, and of similar quantities in other models. 
The methods possess several weaknesses: 
Conditions for validity of the approximations are not known, 
magnitudes of approximation errors are not easily evaluated, 
spurious solutions are generated that require large efforts to eliminate, 
expressions for the approximations are in many cases too complex to be useful. 
We describe an alternative method that provides improvements in these regards. 
The new method leads to asymptotic approximations of the first few cumulants that are explicit in the model's parameters. 
We analyze the univariate stochastic logistic Verhulst model and a bivariate stochastic epidemic SIR model with the new method.  
Errors that were made in early applications of moment closure to the Verhulst model are explained and corrected. 


\end{abstract}

\section{Introduction}

A common situation in the study of many stochastic models is that one wishes to determine 
values of a few low order moments or cumulants of some random variable of interest. 
It is furthermore common that one can derive a system of ordinary differential equations 
(ODEs) for these moments or cumulants. 
(We deal mainly with cumulants in what follows.) 
It is natural to proceed to solve the ODEs for the unknown cumulants. 
This step is, however, made complicated by the fact that the equations are not closed in the interesting cases where the transition rates are nonlinear. 
This means that the equations for the cumulants up to a given order $m$ contain cumulants of order larger than $m$.   
Cumulant closure has then been used to close the equations before they are solved. 
It takes the form for any $m\ge 1$ that all cumulants of order exceeding $m$ that appear 
in the system of ODEs are expressed in terms of the cumulants of order at most equal to $m$. 
A common way of achieving this is to make a distributional assumption about the random variable 
of interest, and another one is to use cumulant neglect and put all cumulants of order exceeding 
$m$ equal to zero. 
By solving the closed system of ODEs for the cumulants up to order $m$ we are led to approximations 
of these cumulants. 
We shall be particularly interested in stationary values of the cumulants. 
They are found as coordinates of critical points of the system of ODEs.

Numerical studies have shown that there are situations where the approximations 
do not work at all, e.g. by leading to negative variances. 
In other cases, however, one has found that cumulant closure leads to quite reasonable 
approximations of the cumulants. 
Several different distributional assumptions are possible in cases where the approximations appear 
acceptable.
The different assumptions lead to different approximations, but no theoretical basis exists for choice among them. 
The magnitudes of the errors caused by cumulant closure can not be evaluated from the method itself. 
In many cases spurious solutions appear that require additional efforts to study and eliminate. 
Expressions for the resulting approximations can be used for numerical evaluations, but are often 
too complex to give insight into the dependence on model parameters.

The main aim of the present paper is to describe an alternative to moment closure that has several properties that make it more attractive than moment closure from a modelling standpoint. 
The new method starts out with the same system of ODEs for the low order cumulants as the 
moment closure method. 
However, it deviates from the moment closure approach by going directly for approximations 
of the cumulants without taking any step for closing the system of ODEs. 
In this way one avoids the undesirable consequencies caused by the ad-hoc nature of the 
moment closure assumptions.  
The cumulant approximations that are produced by the new method are asymptotic as 
some parameter takes large values. 
We describe the new method in detail by analyzing both the well-studied univariate stochastic logistic 
Verhulst population model, and a bivariate stochastic epidemic SIR model. 
We also explain and correct a conceptual error that was committed in early work applying 
moment closure to the Verhulst model.

A broad description of moment closure methods for a host of mathematical models is given 
in the recent review by Kuehn (2016). 
An early description dealing wih the foundations of the method is contained in the 
basic paper by Whittle (1957).  
He uses $m=2$ and achieves cumulant closure by making a normal approximation. 
Since cumulants of order exceeding 2 are equal to zero for a normally distributed random 
variable, this is equivalent to cumulant neglect.

The alternative to moment closure that we describe here will in some cases require a reparametrization of the model so that one can identify a parameter that takes large values. 
This parameter will then serve the important role as the one for which the results derived 
serve as asymptotic approximations for large values of the parameter in question. 
In the case of the stochastic logistic Verhulst model that is treated below, we shall use the 
maximum population size $N$ for this purpose.
In the bivariate SIR model, the large parameter is again denoted by $N$, and interpreted as the 
expected population size. 
The new method needs information about the orders of magnitude of the cumulants 
that appear in the system of ODEs.

The univariate stochastic logistic Verhulst population model is dealt with in Sections 2--4. 
The model formulation is treated in Section 2, where we also describe two different 
parametrizations that have appeared in the literature. 
The more recent one is necessary for application of the alternative to cumulant closure
that we describe here. 
Section 3 deals with two approaches that have been used in the study of the 
stochastic logistic Verhulst model with cumulant closure methods. 
We show that errors were made in some of the early work. 
The application of the new method to the stochastic logistic Verhulst model is treated in Section 4. 
Asymptotic approximations of the first 3 cumulants of the so-called quasi-stationary distribution (QSD) 
are derived. 
This is a stationary distribution of the state variable conditioned on nonextinction. 
It is useful for understanding the behavior of the model whenever extinction has not occured.  
Conditions on the parameter space for validity of the results are given, as well as magnitudes 
of error terms. 
A bivariate stochastic SIR model is studied in Section 5. 
Essentially the same ideas used for the analysis of the univariate Verhulst model are applied. 
Asymptotic approximations of the 5 cumulants of orders 1 and 2 of the bivariate QSD are derived. 
The paper ends by some concluding remarks in Section 6.

\section{The stochastic logistic Verhulst model: Two model formulations}

The stochastic logistic Verhulst model is formulated as a birth-death process $\{X(t),t\ge 0\}$. 
The hypotheses of the model are summarized by descriptions of the population birth-rate 
$\lambda_n$ and the population death-rate $\mu_n$ as functions of the state $n$ of the process. 
Two different formulations of these transition rates have been used in the literature. 
We describe both of them, and the corresponding parameter spaces.

Realistic population models account for the fact that population growth is density 
dependent, in the sense that the net growth rate per individual is a decreasing function 
of the population size. 
The classical deterministic population model formulated by Verhulst (1838) allowed for 
density dependence by hypothesizing the net growth rate per individual to be a linearly 
decreasing function of the population size. 
A host of additional deterministic population models whose growth rates are nonlinearly 
decreasing functions of the population size is studied by Tsoularis and Wallace (2002).

Two different ways of formulating the stochastic logistic Verhulst model have appeared
in the literature. 
In the first one historically, the population birth-rate $\lambda_n$ is expressed as 
\begin{equation} \label{2.1}
   \lambda_n = \begin{cases}
               & (a_1 - b_1n)n,  \quad n=0,1,\dots, [a_1/b_1],  \\
               & 0, \quad n \ge a_1/b_1,
               \end{cases}
\end{equation}
and the population death-rate $\mu_n$ as 
\begin{equation} \label{2.2}
   \mu_n = (a_2 + b_2 n)n, \quad n=0,1,2,\dots.
\end{equation}

As for any model, it is important to describe its state space and its parameter space. 
The state space in this formulation has not always been specified, but it 
appears in most cases to consist of all non-negative integers.  
The parameter space consists of the four parameters $a_1, a_2, b_1, b_2$. 
They are variously assumed to be non-negative or strictly positive.

This formulation was introduced by Bartlett, Gower, and Leslie (1960), 
and it has after this been followed by Matis and Kiffe (1996), 
by Singh and Hespanha (2007), and by Renshaw (2011). 
Krishnarajah, Cook, Marion, and Gibson (2005) study the SIS model, which is 
a special case of the more general logistic Verhulst model, as mentioned below. 
They use a similar model formulation, but with the particular feature that the 
maximum population size $N$ is one of the parameters.

A second formulation of this stochastic model was given by N{\aa}sell 
(2001), (2011), as follows. 
The population birth-rate $\lambda_n$ was expressed as 
\begin{equation} \label{2.3}
   \lambda_n = \mu R_0 \left(1 - \alpha_1 \frac{n}{N}\right)n,  
                 \quad n=0,1,\dots, N-1,  \quad \lambda_N = 0, 
\end{equation}
and the population death-rate $\mu_n$ as
\begin{equation} \label{2.4}
   \mu_n = \mu \left(1+\alpha_2 \frac{n}{N} \right)n, \quad n=0,1,\dots,N.
\end{equation}

In this formulation, $\alpha_1$ was assumed to take values in the closed unit interval: 
$0\le \alpha_1 \le 1$.  
In what follows we shall deviate from this formulation by restricting attention to the case $\alpha_1=1$. 
The reason for this is that if $\alpha_1<1$, then the birth-rate $\lambda_n$ would experience 
a large and biologically unmotivated change when $n$ is increased from $N-1$ to $N$. 
In the model that we study here we shall therefore use the following definitions of 
$\lambda_n$ and $\mu_n$:
\begin{equation} \label{2.5}
    \lambda_n = \mu R_0 \left(1 -\frac{n}{N}\right)n,  \quad n=0,1,\dots, N,  
\end{equation}
and 
\begin{equation} \label{2.6}
   \mu_n = \mu \left(1+\alpha \frac{n}{N} \right)n, \quad n=0,1,\dots,N.
\end{equation}  

The state space of the process in this case is finite, being equal to $\{0,1,\dots,N\}$. 
The parameter space for this model formulation consists of the four parameters 
$N, R_0, \alpha, \mu$. 
Among these, $N$ is a large positive integer that represents the maximum 
population size, $R_0$ and  $\alpha$ are dimensionless parameters, while 
$\mu$ is a positive death rate with the dimension inverse time. 
The parameter $R_0$ is a positive threshold parameter discussed below, 
while $\alpha$ is a non-negative constant.

We note that the SIS model for transmission of infection without immunity 
in a constant population of hosts can be seen as a special case of this formulation of 
the logistic model with $\alpha=0$.  
A consequence of this is that the approximations that we derive for the 
cumulants of the quasi-stationary distribution of the logistic Verhulst model are also 
valid for the SIS model.

The model formulation based on the transition rates in \eqref{2.3}-\eqref{2.4} 
has been studied by the aid of moment closure methods by N{\aa}sell (2003a, 2003b), 
Newman, Ferdy, and Quince (2004), Clancy (2012), and Martins, Pinto, and Stollenwerk (2012), 
where the latter authors restrict themselves to the SIS model.

The parameter $R_0$ serves the important role of identifying a threshold 
at $R_0=1$ for the deterministic version of the model. 
The solution of the deterministic model shows qualitatively different 
behaviors in the two parameter regions above threshold ($R_0>1$) and at or below 
threshold ($0<R_0\le 1$). 
Indeed, given a positive initial value, the population size is predicted by the 
deterministic model to approach a positive level above the threshold, while it 
is predicted to go extinct at or below the threshold.

In similarity to this, the stochastic model shows qualitatively different 
behaviors in three parameter regions, as shown by N{\aa}sell (2001). 
Thus, for $R_0>1$ and large $N$, we find for the stochastic version of the model 
that the time to extinction is exponentially large, and that the quasi-stationary 
distribution (QSD) of the population size is approximately normal in its body. 
In contrast to this we find for $R_0<1$ and large $N$ that the time to  
extinction is short and that the QSD is approximately geometric in its left tail.  
In fact, the time to extinction is then so short that there is not always 
enough time (depending on initial conditions) for the distribution to approach 
quasi-stationarity before extinction occurs. 
This means that the concept of quasi-stationarity is not interesting or useful 
in this parameter region. 
A third parameter region is defined when $R_0$ is close to its threshold value one. 
The asymptotic study that we advocate requires a reparametrization in this region. 
It is achieved by defining a new parameter $\rho$ by the expression  
\begin{equation} \label{2.7}
    \rho = \frac{R_0-1}{\sqrt{1+\alpha}}\sqrt{N},  
\end{equation} 
and keeping $\rho$ fixed as $N$ grows toward infinity. 
This parameter region is referred to as a transition region. 
Clearly, the time to extinction is moderately large in this region, while the 
QSD makes a transition from being approximately normal in its body 
to being approximately geometric in its left tail when $R_0$ is reduced from above
to below the value 1.  
Properties of the QSD are discussed in N{\aa}sell (2001, 2011).

The step from the original model formulation given by \eqref{2.1}-\eqref{2.2}  
to the second one in \eqref{2.5}-\eqref{2.6} involves a reparametrization. 
It makes use of the concepts of dimensional analysis and scaling. 
These ideas are common in physical modelling, but appear to be less 
well-known in stochastic modelling of population processes. 
A discussion for the latter case is given by N{\aa}sell (2002).

\section{The stochastic logistic Verhulst model: Two moment closure approaches}

Two random variables have been used to study the stochastic logistic Verhulst model 
with the aid of moment closure methods. 
One of them uses the state variable $X(t)$ of the process. 
It has an absorbing state at the origin. 
Absorption at this point corresponds to extinction of the population studied. 
The second random variable that has been used, denoted $X^Q(t)$, is defined by 
conditioning $X(t)$ on non-extinction.

Both the moment closure method and its new alternative work with a system of ODEs 
for the first few cumulants of the random variable that is studied. 
Goals of both methods are to derive approximations of the first few cumulants of the 
stationary distribution of the random variable that one is concerned with.  
Cumulants of the stationary distribution appear as coordinates of a critical point of the system
of ODEs.

The stationary distribution of $X(t)$ is degenerate with probability one at the origin. 
All its cumulants are therefore equal to zero. 
It is therefore not necessary to use the associated system of ODEs for the cumulants at all 
to determine these cumulants. 
However, early applications of moment closure to the stochastic logistic Verhulst model erred by 
not observing this simple fact. 
Instead they studied the system of ODEs for the unconditioned random variable $X(t)$. 
In contrast to this, the stationary distribution of the conditioned random variable 
$X^Q(t)$ equals the QSD, and is of substantial interest.

We find in particular that the stationary values of the cumulants of the conditioned random 
variable $X^Q(t)$ are equal to the cumulants of the QSD.  
In this case it is entirely appropriate to use the system of ODEs for determining 
the first few cumulants of the QSD. 
The alternative method that we advocate here works only with the conditioned 
random variable $X^Q(t)$.

It is quite remarkable that cumulant closure in early approaches using 
$X(t)$ leads to acceptable approximations of the cumulants of the 
stationary distribution of the conditioned random variable $X^Q(t)$. 
An explanation for this fact is given in the next section of the paper.

The papers by Bartlett, Gower, and Leslie (1960), Matis and Kiffe (1996), 
N{\aa}sell (2003a), Newman, Ferdy, and Quince (2004), 
Krishnarajah, Cook, Marion, and Gibson (2005), and the book by  Renshaw (2011) 
all work with the random variable $X(t)$, while N{\aa}sell (2003b), 
Singh and Hespanha (2007), Clancy (2012), and Martins, Pinto, and 
Stollenwerk (2012) use the conditioned random variable $X^Q(t)$.

The new method presented here deals only with the conditioned random variable 
$X^Q(t)$, and it uses only the second of the two parameter spaces described above. 
The important difference from the classical moment closure method is that 
the steps taken to achieve moment closure are replaced by assumptions about the 
form of the asymptotic approximations of the first few cumulants of the QSD for large $N$.

Cumulant closure applied to the two random variables $X(t)$ and 
$X^Q(t)$ was studied by N{\aa}sell in (2003a) and (2003b), respectively.  
The fact that the first one of these studies was completely unnecessary had not 
been recognized at that time. 
In both cases we used the second of the two parametrizations of Section 2, 
and in both cases we used a two-step process: 
Cumulant closure was followed by asymptotic approximation. 
It was then argued that the final result provided asymptotic 
approximations of the first few cumulants studied. 
However, the logic behind this conclusion can be questioned. 
The reason for this is that no knowledge is available about the magnitudes 
of the errors associated with the moment closure method. 
It does not make sense to study asymptotic approximations of these results, 
since they involve approximation errors of unknown magnitudes. 
In contrast to this, we claim that the new method that derives asymptotic 
approximations of the original cumulants without any intermediate step of 
moment closure gives results in which the magnitudes of the error terms are 
known. 
The search for moment closure is replaced by a search for asymptotic approximation.

Matis and Kiffe (1996) have given differential equations for the first three cumulants of 
the unconditioned random variable $X(t)$. 
We express their results with the aid of the second parametrization of Section 2. 
The results can be written as follows:

\begin{align} \label{3.1}
    \kappa_1'(t) & = \mu A(t), \\ \label{3.2}
    \kappa_2'(t) & = \mu B(t), \\  \label{3.3}
    \kappa_3'(t) & = \mu C(t), 
\end{align}

{\allowdisplaybreaks  

where
\begin{align} \label{3.4}
   A(t) & =  (R_0-1) \kappa_1(t) -\frac{R_0 + \alpha}{N} 
     [\kappa_1^2(t)+ \kappa_2(t)], 
\\ \label{3.5}
   B(t) & =  (R_0+1)\kappa_1(t) + 2 (R_0-1) \kappa_2(t) 
      - \frac{R_0 - \alpha}{N} [\kappa_1^2(t) + \kappa_2(t)]   \\* \notag
      &  \phantom{abc} - \frac{R_0 + \alpha}{N} 
       [4 \kappa_1(t) \kappa_2(t) + 2 \kappa_3(t)],  
\\ \label{3.6}
  C(t) & =  (R_0-1) [\kappa_1(t) + 3 \kappa_3(t)] + 3 (R_0+1) \kappa_2(t) \\* \notag
     & \phantom{abc}  - \frac{R_0 - \alpha}{N} [6 \kappa_1(t) \kappa_2(t) + 3 \kappa_3(t)]
        \\* \notag 
      &   \phantom{abc} - \frac{R_0 + \alpha}{N} [\kappa_1^2(t) 
            + 6 \kappa_1(t) \kappa_3(t) + \kappa_2(t) + 6 \kappa_2^2(t) + 3 \kappa_4(t)].
\end{align}
} 
 
We  note that these equations are not closed, since $C(t)$ depends on $\kappa_4(t)$. 
One way of achieving cumulant closure is to use cumulant neglect and put $\kappa_4(t)=0$.  
A critical point of the resulting system of equations for the first 3 cumulants is found 
by putting $\kappa_4(t)=0$ and solving the equations $A(t)=B(t)=C(t)=0$ for the 
stationary values of the first 3 cumulants $\kappa_1, \kappa_2, \kappa_3$.  
There are two solutions of this mathematical problem. 
One is simply $\kappa_1=\kappa_2=\kappa_3=0$, while the second one has $\kappa_1>0$. 
The first of these solutions corresponds to the known stationary distribution of the 
random variable $X(t)$, namely the degenerate one with probability one at the origin, while 
the second one leads to several critical points of the system \eqref{3.1}-\eqref{3.3}. 
They all correspond to spurious solutions that have no correspondence to stationary 
distributions of the nonconditioned random variable $X(t)$. 
Early studies of cumulants for the stochastic logistic model made 3 mistakes. 
The first one was to consider the stationary cumulants of the random variable $X(t)$ instead of
$X^Q(t)$. 
The second one was to use the differential equations \eqref{3.1}--\eqref{3.3} instead of 
immediately noting that the stationary cumulants of $X(t)$ are all equal to zero, 
and the third one was to search for solutions of the mathematical problem described 
above with positive mean $\kappa_1$.

The right way to go is to determine stationary cumulants of $X^Q(t)$, which coincide with 
cumulants of the QSD.  
Differential equations for the first three cumulants of the conditioned random variable 
$X^Q(t)$ have been derived by N{\aa}sell (2003b).  
The resulting system of ODEs for the first three cumulants is similar to 
the system given in \eqref{3.1}-\eqref{3.3}. 
It takes the form

\begin{align}  \label{3.7}
   \kappa_1'(t) & = \mu A(t) + \mu_1 q_1(t) \kappa_1(t), \\ \label{3.8}
   \kappa_2'(t) & = \mu B(t) +  \mu_1 q_1(t) [\kappa_2(t) - \kappa_1^2(t)], \\ \label{3.9}
   \kappa_3'(t) & = \mu C(t)  + \mu_1 q_1(t) [\kappa_1^3(t) - 3\kappa_1(t)\kappa_2(t) + \kappa_3(t)].
\end{align}

In the last terms in the right-hand sides of the 3 equations \eqref{3.7}--\eqref{3.9}, 
we find from the definition of $\mu_n$ in \eqref{2.4} that $\mu_1 = \mu(1+\alpha/N)$.
Furthermore, $q_1(t)$ is used to denote the probability that the conditioned 
random variable $X^{Q}(t)$ takes the value one. 
The steady-state value of this probability is shown by N{\aa}sell (2001) to be exponentially
small in the parameter region where $R_0>1$. 
The corresponding terms can therefore be ignored when we search for asymptotic approximations 
of the stationary cumulants.

Cumulant closure of the system of differential equations \eqref{3.7}-\eqref{3.9} is as above 
achieved by cumulant neglect. 
Thus, we find that a critical point of the resulting system of equations for the three 
cumulants of $X^Q(t)$ is found by solving the equations $A(t)=B(t)=C(t)=0$, with 
$\kappa_4(t)=0$. 
In this way we are led to the same mathematical problem as above where we were concerned with 
the cumulants of $X(t)$, with the exception that in this case it is correct to search for a solution with  $\kappa_1>0$.

The point where $\kappa_1=\kappa_2=\kappa_3=0$ is excluded in this case, since it does not 
correspond to any stationary distribution of the conditioned random variable $X^Q(t)$. 
On the other hand, as is further discussed by N{\aa}sell (2003b), additional critical points are 
associated with any solution where $\kappa_1 \ne 0$. 
They are determined from the roots of an equation of fourth degree. 
Further study of the roots of this equation is required in order to identify one of these roots 
as acceptable, while the remaining three of them are associated with spurious solutions. 
The investigations for this require determinations of stabilities of the corresponding critical points.

We note that the two tasks of determining stationary values of the first three cumulants of $X(t)$ 
and of $X^Q(t)$ satisfy similar mathematical problems. 
The only exception is that one should choose the critical point with $\kappa_1=0$ in the first case, 
and the one with $\kappa_1>0$ in the second case. 
Since the incorrect choice of taking $\kappa_1>0$ was made in the first case, we find that the two mathematical problems are identical.  
This allows us to understand the surprising fact that the cumulant approximations that
were derived by applying the cumulant closure method to $X(t)$ agree with the approximations 
derived with the same method applied to $X^Q(t)$, with $R_0>1$.

\section{The stochastic logistic Verhulst model: Asymptotic approximations}

This section is used to establish asymptotic approximations of the first 
three cumulants of the quasi-stationary distribution of the stochastic logistic Verhulst model 
in the parameter region where $R_0>1$. 
We base our results on assumptions concerning the forms of asymptotic approximations 
of the first four cumulants for large values of $N$. 
Before such assumptions can be established we need information about the orders of 
magnitude of these cumulants. 
One basis for these assumptions is the numerical evaluations given in Table 1. 
The table shows the values of the first 4 cumulants of the 
quasi-stationary distribution in each of the 3 parameter regions, and
for 3 different values of $N$. 
The results in the transition region ($R_0=1$) and in the region distinctly 
below threshold ($R_0=0.4$) will not be used here, but are shown for their 
independent interest. 
The evaluations in the table have been made for the SIS model, 
with $\alpha=0$.

\begin{table}[h]
  \begin{center}
    \begin{tabular}{ | c | c | r | r | r | }
       \hline
    $R_0$&  Cumulant  &N=100& N=200 &N=400\\ \hline
     0.4 & $\kappa_1$ &1.64 & 1.65  & 1.66 \\
     0.4 & $\kappa_2$ &1.02 & 1.06  & 1.09 \\
     0.4 & $\kappa_3$ &2.22 & 2.39  & 2.49 \\ 
     0.4 & $\kappa_4$ &6.64 & 7.49  & 7.98 \\ \hline
     1   & $\kappa_1$ &7.03 & 9.80  & 13.7 \\
     1   & $\kappa_2$ &27.3 & 55.9  & 114  \\
     1   & $\kappa_3$ & 160 & 476   & 1394 \\ 
     1   & $\kappa_4$ & 899 & 3983  & 17072 \\ \hline
     2   & $\kappa_1$ &48.9 & 99.0  & 199 \\
     2   & $\kappa_2$ &52.3 & 102   & 202 \\
     2   & $\kappa_3$ &-58.2& -107  & -206 \\ 
     2   & $\kappa_4$ & 95.1& 133   & 229 \\ \hline
    \end{tabular}
    \vskip 4mm
    \caption{Numerical evaluations of the first 4 cumulants of the 
     QSD of the SIS model.
     Results are shown in each of the 3 parameter regions, and for 3 
     different $N$-values. 
     }
  \end{center}
\end{table}

From the entries in the table we see that the first 4 cumulants are practically 
independent of $N$ for $R_0=0.4$. 
This indicates that the cumulants are of order 1 when $R_0<1$. 
In contrast to this, we see that the first 4 cumulants are multiplied by 
approximately 2 when $N$ is doubled and $R_0=2$.  
We interpret this as an indication that the first four cumulants are of order $N$ when $R_0>1$.  
A rather different behavior of the cumulants is seen when $R_0=1$. 
We note here that a doubling of $N$ has the consequence that $\kappa_1$ is 
multiplied by approximately $\sqrt{2}$, $\kappa_2$ by 2, $\kappa_3$ by 
$2\sqrt{2}$, and $\kappa_4$ by 4. 
This indicates that the $i$th cumulant in the transition region is of the order of $N^{i/2}$, 
at least for the $i$-values 1, 2, 3, and 4.

We proceed to study the first 3 cumulants of the QSD in the parameter region where 
$R_0>1$ and $N$ is large. 
As a start we need information about the orders of magnitude of the first 4 cumulants. 
We assume that they are all of the order of $N$. 
This assumption is clearly supported by the numerical results discussed above. 
Additional support for this assumption is given by N{\aa}sell (2001), where it is shown that the QSD 
for $R_0>1$ is approximately normal in its body, with mean equal to $(R_0-1)N/(R_0+\alpha)$ and variance equal to $(1+\alpha)R_0 N/(R_0+\alpha)^2$. 
These approximations clearly support the assumptions that the first cumulant $\kappa_1$, 
which equals the mean, and the second cumulant $\kappa_2$, which equals 
the variance, are of the order of $N$. 
The assumptions that the third and fourth cumulants are also of order $N$ are equivalent to the  assumptions that neither of them is of larger order than $N$. 
Any such larger order would be completely at odds with the finding that the QSD is approximately 
normal in its body, since all cumulants of order exceeding 2 are equal to 0 for the 
normal distribution.

The next step in the study of the first 3 cumulants is to introduce assumptions concerning 
the first few terms of the first 4 cumulants of the QSD. 
These assumptions take the following form: 
{\allowdisplaybreaks
\begin{align}  \label{4.1}
  \kappa_1 & = x_1 N + x_2 + \frac{x_3}{N} + \text{O}\left(\frac{1}{N^2}\right), 
     \quad R_0>1, \\ \label{4.2}
  \kappa_2 & = x_4 N + x_5 + \text{O}\left(\frac{1}{N}\right), \quad R_0>1, 
    \\ \label{4.3}
  \kappa_3 & = x_6 N + \text{O}(1), \quad R_0>1, \\ \label{4.4}
  \kappa_4 & = x_7 N + \text{O}(1), \quad R_0>1. 
\end{align}

In these expressions, we have introduced the 7 quantities $x_1$--$x_7$. 
It is important that they all are independent of $N$. 
Knowledge about the first 6 of them gives asymptotic approximations of the first 3 cumulants.

One may be tempted to include equal numbers of terms in the assumed forms for all cumulants. 
However, the next step in our development is to determine asymptotic 
approximations of the quantities $A$, $B$, and $C$. 
By accounting for orders of magnitude of the errors committed to each of these
quantities by the assumed forms, one would find that e.g. a third term 
in the expression for $\kappa_2$ would contribute a term to $A$ that has 
the same magnitude as the error term already contributed to $A$ by using the 
three terms assumed for $\kappa_1$.

It is straightforward to insert the assumed asymptotic expressions for the 
first 4 cumulants into the expressions \eqref{3.4}--\eqref{3.6} for the 
quantities $A$, $B$, and $C$, and to derive the resulting asymptotic 
expressions for $A$, $B$, and $C$. 
The results take the following form:

\begin{align} \label{4.5}
   A & = A_1 N + A_2 + \frac{A_3}{N} + \text{O}\left(\frac{1}{N^2} \right), 
     \quad R_0>1,  \\ \label{4.6}
   B & = B_1 N + B_2 + \text{O}\left(\frac{1}{N} \right),   \quad R_0>1,  \\ \label{4.7}
   C & = C_1 N + \text{O}(1),  \quad R_0>1, 
\end{align}
where
\begin{align} \label{4.8}
   A_1 & = (R_0-1) x_1 -(R_0+ \alpha) x_1^2, \\ \label{4.9}
   A_2 & = (R_0-1) x_2 - (R_0 + \alpha)(2x_1 x_2 + x_4),   
      \\ \label{4.10}
   A_3 & = (R_0-1) x_3 - (R_0 + \alpha) (2 x_1 x_3 +x_2^2 + x_5), 
      \\ \label{4.11}
   B_1 & = (R_0+1) x_1 + 2(R_0-1) x_4 -(R_0 - \alpha)x_1^2 
                - 4(R_0 + \alpha)x_1 x_4, \\ \label{4.12}
   B_2 & = (R_0+1) x_2 + 2(R_0-1) x_5 - (R_0 - \alpha) 
       (2 x_1 x_2 + x_4) \\* \notag
   & \phantom{abc}- (R_0 + \alpha) (4 x_1 x_5 + 4 x_2 x_4 + 2 x_6),  
      \\ \label{4.13}
   C_1 & = (R_0-1) (x_1 + 3 x_6) + 3 (R_0+1) x_4 
       - 6 (R_0 - \alpha) x_1 x_4 \\* \notag
     & \phantom{abc} - (R_0 + \alpha) (x_1^2 + 6 x_1 x_6 + 6 x_4^2).
\end{align}

These expressions can be used to form 6 equations by setting each of the
expressions equal to zero. 
We note also that the 6 coefficients $x_1$--$x_6$ will determine the first term or 
terms that give asymptotic appproximations for the first 3 cumulants. 
Thus, since the number of equations equals the number of unknown quantities, solution appears possible. 
By further inspection we find that the values of the $x_i$ can be determined sequentially, as follows:

First, the equation $A_1=0$ is solved for $x_1$. 
This equation has two roots, namely $x_1=0$ and $x_1=(R_0-1)/(R_0+\alpha)$. 
Among them, we exclude $x_1 = 0$ as the only spurious solution that appears in 
this method. 
After this, $B_1=0$ is solved for $x_4$, and $A_2=0$ is solved for $x_2$. 
The result so far determines the critical point of the system of ODEs for the 
first two cumulants. 
It gives a two-term approximation of $\kappa_1$ and a one-term approximation 
of $\kappa_2$. 
To work with the critical point of the system of ODEs for the first three 
cumulants, we continue and
solve $C_1=0$ for $x_6$, $B_2=0$ for $x_5$, and $A_3=0$ for $x_3$. 
We note that all equations after the first one are linear and elementary to solve. 
The resulting values of the coefficients $x_i$ are as follows:
\begin{align} 
   x_1 & = \frac{R_0-1}{R_0+\alpha}, \\
   x_2 & = - \frac{(1+\alpha)R_0}{(R_0+\alpha)(R_0-1)}, \\
   x_3 & = -\frac{(1+\alpha)(R_0+1)R_0}{(R_0-1)^3}, \\
   x_4 & = \frac{(1+\alpha)R_0}{(R_0+\alpha)^2}, \\ 
   x_5 & = \frac{(1+\alpha)((R_0^2+\alpha)R_0}{(R_0+\alpha)^2(R_0-1)^2}, \\ 
   x_6 & = -\frac{(1+\alpha)(R_0-\alpha)R_0}{(R_0+\alpha)^3}. 
\end{align}

The resulting values of $x_1$--$x_6$ inserted into 
\eqref{4.1}--\eqref{4.3} lead to the following asymptotic approximations of 
the first 3 cumulants in the parameter region where $R_0>1$:

\begin{align} \label{4.20}
   \kappa_1 & = \frac{R_0-1}{R_0 + \alpha} N 
     -\frac{(1 + \alpha) R_0}{(R_0 + \alpha)(R_0-1)}  
     - \frac{(1 + \alpha) (R_0+1) R_0}{(R_0-1)^3}\frac{1}{N} \\* \notag
      & \phantom{abc}+ \text{O}\left( \frac{1}{N^2}\right),   \quad R_0>1,  
         \\ \label{4.21}
  \kappa_2 & =  \frac{(1 + \alpha) R_0}{(R_0 + \alpha)^2} N 
       + \frac{(1+\alpha)(R_0^2+\alpha) R_0} {(R_0 + \alpha)^2 (R_0-1)^2}  
         + \text{O}\left(\frac{1}{N}\right),  \quad R_0>1,  
        \\ \label{4.22}
   \kappa_3 & =  - \frac{(1 + \alpha) (R_0 - \alpha) R_0}
        {(R_0 + \alpha)^3} N + \text{O}(1), \quad R_0>1. 
\end{align}
} 

Maple has been used to aid in the book-keeping required to derive these results. 
These approximating expressions for the first 3 cumulants of the QSD are not new. 
They agree formally with results derived using moment closure followed by asymptotic 
approximation and given by N{\aa}sell (2003b). 
The new method introduced here has thus been used to show that the early results based on 
moment closure have the desirable property of being asymptotic.

We note that the two competing methods both require a bound on $\kappa_4$. 
The new method allows the absolute value of $\kappa_4$ to grow with $N$, but not faster than
$\kappa_4=\text{O}(N)$.
This represents an appreciable relaxation of the cumulant closure requirement that 
$\kappa_4(t)=0$.

Numerical evaluations of the error terms in the approximations \eqref{4.20}--\eqref{4.22} 
are given in Table 2. 
They are consistent with the results that the error terms of $\kappa_i$ are of the orders of 
$1/N^{3-i}$ for the $i$-values 1, 2, and 3. 
To see this, we note that the error terms for $\kappa_1, \kappa_2, \kappa_3$ 
are divided by approximately 4, 2, and 1, respectively, when $N$ is doubled.

\begin{table}[h]
  \begin{center}
    \begin{tabular}{ | c | c | r | r | r | }
       \hline
    $R_0$&  Cumulant  &N=100& N=200 &N=400\\ \hline
     2   & $\kappa_1$  &-0.0095 & -0.0021 & -0.00050 \\
     2   & $\kappa_2$  &0.33      & 0.15     &  0.069 \\
     2   & $\kappa_3$  &-8.2      & -6.9      & -6.4 \\ \hline
    \end{tabular}
    \vskip 4mm
    \caption{Numerical evaluations of the error terms of the approximations 
     \eqref{4.20}--\eqref{4.22} for the first 3 cumulants of the 
     QSD of the SIS model. 
     Results are shown for $R_0=2$, and for 3 different $N$-values. 
     }
  \end{center}
\end{table}

If an additional term in the asymptotic approximation of one of the cumulants 
$\kappa_1-\kappa_3$ is required, then it is necessary to increase the number 
of cumulant ODEs that are analyzed from 3 to 4. 
The same obviously holds if one wishes to derive an approximation of the 
cumulant $\kappa_4$.

Some authors work with the raw moments  $\bar{\mu}_i$ (moments about zero) 
instead of the cumulants $\kappa_i$.  
By using the relations 
\begin{align} \label{4.23}
   \bar{\mu}_1 & = \kappa_1,  \\ \label{4.24}
   \bar{\mu}_2 & = \kappa_2 + \kappa_1^2, \\ \label{4.25}
   \bar{\mu}_3 & = \kappa_3 + 3\kappa_1 \kappa_2 + \kappa_1^3,
\end{align}
one finds from the results in \eqref{4.20}--\eqref{4.22} that 
approximations of the first 3 raw moments of the QSD 
of the stochastic logistic Verhulst model can be written as follows:
{\allowdisplaybreaks
\begin{align} \label{4.26}
  \bar{\mu}_1 & = \frac{R_0-1}{R_0 + \alpha} N
           -\frac{(1 + \alpha) R_0}{(R_0 + \alpha)(R_0-1)} 
    \\* \notag
     & \phantom{abc} -\frac{(1 + \alpha) (R_0+1) R_0}{(R_0-1)^3} 
     \frac{1}{N} + \text{O}\left(\frac{1}{N^2}\right),    \quad R_0>1,  
     \\ \label{4.27}
  \bar{\mu}_2 & = \frac{(R_0-1)^2}{(R_0 + \alpha)^2} N^2 
         - \frac{(1 + \alpha) R_0}{(R_0 + \alpha)^2} N   
      \\* \notag
         & \phantom{abc} - \frac{(1 + \alpha) (R_0+1)R_0}
         {(R_0 + \alpha) (R_0-1)^2}
           + \text{O}\left(\frac{1}{N}\right),    \quad R_0>1,  \\  \label{4.28}
  \bar{\mu}_3 & = \frac{(R_0-1)^3}{(R_0 + \alpha)^3} N^3  \\* \notag
             & \phantom{abc} - \frac{(1+\alpha)[R_0^2 + 2(1+\alpha)R_0+\alpha] R_0}
         {(R_0 + \alpha)^3 (R_0-1)} N\\* \notag
             & \phantom{abc}  + \text{O}(1),    \quad R_0>1.
\end{align}
}

Maple has been used to derive these results. 

\begin{table}[h]
  \begin{center}
    \begin{tabular}{ | c | c | r | r | r | }
       \hline
    $R_0$&  Moment  &N=100& N=200 &N=400\\ \hline
     2   & $\bar{\mu}_1$  &-0.0095 & -0.0021 & -0.00050 \\
     2   & $\bar{\mu}_2$  & -0.48     & -0.21    &  -0.10 \\
     2   & $\bar{\mu}_3$  & -27      & -24     & -23 \\ \hline
    \end{tabular}
    \vskip 4mm
    \caption{Numerical evaluations of the error terms of the approximations 
     \eqref{4.26}--\eqref{4.28} for the first 3 raw moments of the 
     QSD of the SIS model. 
     Results are shown for $R_0=2$, and for 3 different $N$-values. 
     }
  \end{center}
\end{table}

Numerical evaluations of the error terms in these expressions are given in 
Table 3. 
They are consistent with the results in \eqref{4.26}--\eqref{4.28} that the 
error terms of $\bar{\mu}_i$ are of the orders of $1/N^{3-i}$ for the 
$i$-values 1, 2, and 3. 
To see this, we note from the table that the error terms for 
$\bar{\mu}_1, \bar{\mu}_2, \bar{\mu}_3$ are divided by approximately 
4, 2, and 1, respectively, when $N$ is doubled.

\section{A bivariate stochastic SIR model}

No new ideas are needed to apply the method described above to bivariate or multivariate models. 
We illustrate this by giving a brief analysis of a bivariate model. 
We choose to deal wth a so-called SIR model with demography. 
This is a model for the transmission of an infection that causes immunity in a population whose size 
is determined by an immigration-death process. 
For historical reasons it has been referred to as the Martini Model. 
The deterministic version of the model has been dealt with by a number of authors, including 
Martini (1921), Lotka (1923, 1956), Hethcote (1974, 1976), and Anderson and May (1991). 
The stochastic version has  been studied by Schenzle (1984), Keeling and Grenfell (1997), 
van Herwaarden and Grasman (1995), and N{\aa}sell (1999, 2005). 
 
\begin{table}[h]
  \begin{center}
    \begin{tabular}{ | c | c | c | }    \hline
    Event & Transition & Transition rate \\ \hline 
    Immigration of susceptible    & $(s,i) \to (s+1,i)$     & $\mu N$       \\ 
    Death of susceptible                & $(s,i) \to (s-1,i)$     & $\mu s$       \\ 
    Infection of susceptible        & $(s,i) \to (s-1,i+1)$  & $\beta si/N$ \\
    Death or recovery of infected & $(s.i) \to (s,i-1)$      & $(\mu+\gamma)i$ \\ \hline
    \end{tabular}
    \vskip 4mm
    \caption{Transition rates for the bivariate SIR model analyzed in Section 5} 
  \end{center}
\end{table}

The model in this section is a bivariate Markov chain $\{(S(t),I(t)), t\ge 0\}$, where $S(t)$ is 
interpreted as the number of susceptible individuals, and $I(t)$ stands for the number of infected individuals, with discrete state space $\{(s,i): s=0,1,2,\dots,i=0,1,\dots\}$ and continuous time.  
It is based on the transition rates given in Table 4. 
The model involves 4 parameters, namely the expected population size $N$, the death rate per individual $\mu$, the contact rate $\beta$, and the recovery rate per infected individual $\gamma$. 
Among these, $N$ is a large positive integer, while the rates $\mu, \beta, \gamma$ are positive 
rates with dimension inverse time. 
A reparametrization is derived in N{\aa}sell (1999, 2005). 
It leads to 2 new parameters $R_0$ and $\alpha$, defined by 
\begin{equation}
   R_0 = \frac{\beta}{\gamma+\mu}
\end{equation}
and 
\begin{equation}
  \alpha = \frac{\gamma+\mu}{\mu}.
\end{equation}
Both of them are positive and dimensionless. 
We note that $\alpha$ is large for the common childhood infections that the model has been 
used to analyze.

A first step in the analysis of this model is given by the diffusion approximation derived by N{\aa}sell (1999). 
It shows that the QSD of the process is approximated by a bivariate normal distribution when $R_0>1$
and $N$ is large. 
This approximating distribution is described by 5 quantities, namely the expectations and variances 
of $S$ and $I$ and their covariance. 
The 2 expectations are equal to the first-order cumulants $\kappa_{10}$ and $\kappa_{01}$, 
while the variances and the covariance are equal to the second-order cumulants 
$\kappa_{20}$, $\kappa_{02}$, and $\kappa_{11}$.  
We quote the following approximations of these 5 cumulants from N{\aa}sell (1999): 
\begin{align} \label{5.3} 
      & \kappa_{10}^{(A)} = \frac{1}{R_0}N, \quad R_0>1, \\ \label{5.4} 
      & \kappa_{01}^{(A)} = \frac{R_0-1}{\alpha R_0}N, \quad R_0>1, \\ \label{5.5} 
      & \kappa_{20} ^{(A)}= \frac{R_0 + \alpha}{R_0^2}N, \quad R_0>1, \\ \label{5.6} 
      & \kappa_{11} ^{(A)}= -\frac{1}{R_0}N, \quad R_0>1, \\ \label{5.7} 
      & \kappa_{02}^{(A)} = \frac{R_0^2 + \alpha(R_0-1)}{\alpha R_0^2}N, \quad R_0>1.
\end{align}

The superscript $(A)$ is used here to indicate that these expressions are approximations of the corresponding cumulants. 
We introduce the assumption that all cumulants of order up to 3 are of order $N$. 
The above results concerning diffusion approximations give strong support to this assumption for the cumulants of order 1 and 2. 
Our assumption for the cumulants of order 3 is equivalent to the assumption that the cumulants of order 3 are not of higher order than $N$. 
Any other behavior would be inconsistent with the facts that the approximating distribution is 
normal, and that all cumulants of the normal distribution of order exceeding 2 are equal to 0.

With these preparatiions we proceed to determine asymptotic approximations of the 
cumulants of the QSD of order 1 and 2. 
The first step is to derive a partial differential equation (PDE) for the cumulant generating function 
$K^{(Q)}(\theta_1,\theta_2,t)$. 
The superscript $(Q)$ is used to indicate that we are conditioning on non-extinction. 
The stationary distribution of this conditioned random variable is equal to the QSD.  
From N{\aa}sell (2005) we find that the PDE for $K^{(Q)}$ can be written as follows:
\begin{multline}  \label{5.8}
   \frac{\partial K^{(Q)}}{\partial t} = \mu 
     \left[\frac{\alpha R_0}{N} (\exp(\theta_2-\theta_1)-1) 
      \left( \frac{\partial^2 K^{(Q)}}{\partial \theta_1 \theta_2} 
     + \frac{\partial K^{(Q)}}{\partial \theta_1} \frac{\partial K^{(Q)}}{\partial \theta_2} \right) \right. \\
     \left.  - (1-\exp(-\theta_1)) \frac{\partial  K^{(Q)}}{\partial \theta_1}
     -\alpha(1-\exp(-\theta_2)) \frac{\partial K^{(Q)}}{\partial \theta_2} 
      +N(\exp(\theta_1)-1) \right] \\
     + \mu \alpha q_{\cdot 1} \left[ 1 - \exp(-K^{(Q)}) \sum_{s=0}^{\infty} q_S(s|1) \exp(s\theta_1)
         \right],
\end{multline}
 where $q_S(s|1) = q_{s1}/q_{\cdot 1}$ is the conditional probability that $S$ takes the value $s$, 
given that $I=1$.   

The cumulant generating function can be expanded in terms of the cumulants.
By spelling out the terms that involve the cumulants of the first 3 orders we get 
\begin{multline}
   K^{(Q)}(\theta_1,\theta_2,t) = \kappa_{10} \theta_1 +\kappa_{01} \theta_2 + \frac12 \kappa_{20}
    \theta_1^2 + \kappa_{11} \theta_1 \theta_2 + \frac12  \kappa_{02} \theta_2^2 \\
    + \frac16 \kappa_{30} \theta_1^3 + \frac12 \kappa_{21} \theta_1^2 \theta_2 
    + \frac12 \kappa_{12} \theta_1 \theta_2^2 + \frac16 \kappa_{03} \theta_2^3 + \dots .
\end{multline}

By using this expansion in the PDE \eqref{5.8} one can derive ODEs for the cumulants of orders 1 and 2. 
This step is more demanding than the corresponding step in the univariate case. 
Maple has been used to help in the bookkeeping necessary for the derivation of the ODEs of these 5 
cumulants. 
The results are given by N{\aa}sell (2005). They can be expressed as follows: 
\begin{align} \label{5.10}
   \kappa_{10}'(t) & = \mu A(t) + \epsilon_A(t), \\ \label{5.11}
   \kappa_{01}'(t) & = \mu B(t) + \epsilon_B(t), \\   \label{5.12}  
   \kappa_{20}'(t) & = \mu C(t) + \epsilon_C(t), \\    \label{5.13}
   \kappa_{11}'(t) & = \mu D(t) + \epsilon_D(t), \\   \label{5.14}
   \kappa_{02}'(t) & = \mu E(t) + \epsilon_E(t).
\end{align}
Here, the functions $A, B, C, D, E$ are as follows:
\begin{align}  \label{5.15}
   A(t) & = N - \frac{\alpha R_0}{N} K_1(t) - \kappa_{10}(t), \\ \label{5.16}
   B(t) & = \frac{\alpha R_0}{N} K_1(t) - \alpha \kappa_{01}(t),  \\ \label{5.17}
   C(t) & = N + \kappa_{10}(t) + \frac{\alpha R_0}{N} (K_1(t) - 2 K_2(t)) - 2 \kappa_{20}(t), \\   \label{5.18}
   D(t) & = \frac{\alpha R_0}{N}(K_2(t)-K_1(t)-K_3(t)) -  (\alpha+1) \kappa_{11}(t), \\ \label{5.19}
   E(t) & = \alpha \kappa_{01}(t) + \frac{\alpha R_0}{N} (K_1(t) + 2 K_3(t)) - 2 \alpha \kappa_{02}(t),
\end{align}
where $K_1, K_2, K_3$ are used to denote the following functions of the cumulants:
\begin{align} \label{5.20}
   K_1(t) & = \kappa_{10}(t) \kappa_{01}(t) + \kappa_{11}(t), \\ \label{5.21}
   K_2(t) & = \kappa_{10}(t) \kappa_{11}(t) + \kappa_{01}(t) \kappa_{20}(t) + \kappa_{21}(t), \\ \label{5.22}
   K_3(t) & = \kappa_{10}(t) \kappa_{02}(t) + \kappa_{01}(t) \kappa_{11}(t) +\kappa_{12}(t). 
\end{align}

The functions $\epsilon_A, \epsilon_B, \epsilon_C, \epsilon_D, \epsilon_E$ appearing in the right-hand sides of equations \eqref{5.10}--\eqref{5.14} are all  proportional to the probability $q_{\cdot 1}(t)$, 
The stationary value of this probability is exponentially small for $R_0>1$ and large $N$. 
Stationary values of the 5 cumulants $\kappa_{10}, \kappa_{01}, \kappa_{20}, \kappa_{11}, 
\kappa_{02}$ are found as coordinates of a critical point of the ODEs \eqref{5.10}--\eqref{5.14}. 
In determining asymptotic approximations of these critical points for $R_0>1$ we can therefore 
ignore the second terms in each of the right-hand sides of the equations \eqref{5.10} -- \eqref{5.14}.

It is seen from \eqref{5.15}--\eqref{5.22} that the functions $A, B, C, D, E$ involve the 2 cumulants 
$\kappa_{10}$ and $\kappa_{01}$ of order 1, the 3 cumulants $\kappa_{20}, \kappa_{11}$, and  $\kappa_{02}$ of order 2, and also the 2 cumulants $\kappa_{21}$ and $\kappa_{12}$ of order 3. 
In order to proceed, we make use of the assumption that the cumulants of the first 3 orders are all 
of the order of $N$ for $R_0>1$. 
We can then introduce the following assumptions concerning these cumulants for the QSD in the parameter region where $R_0>1$ and $N$ is large: 
\begin{align}
    \kappa_{10} & = x_1 N + x_2 + \text{O}\left(\frac{1}{N}\right), \quad R_0>1, \\
    \kappa_{01} & = x_3 N + x_4 + \text{O}\left(\frac{1}{N}\right), \quad R_0>1, \\
    \kappa_{20} & = x_5 N + \text{O}(1), \quad R_0>1, \\
    \kappa_{11} & = x_6 N + \text{O}(1),  \quad R_0>1,  \\
    \kappa_{02} & = x_7 N + \text{O}(1), \quad R_0>1, \\
    \kappa_{21} & = x_8 N + \text{O}(1), \quad R_0>1, \\
    \kappa_{12} & = x_9 N + \text{O}(1), \quad R_0>1. 
 \end{align}

It is straightforward to insert these asymptotic expressions into the expressions \eqref{5.15}--\eqref{5.19} for the functions $A, B, C, D, E$, and derive resulting asymptotic 
approximations of them. 
The results can be written
\begin{align}
    A & = A_1 N + A_2 + \text{O}\left(\frac{1}{N}\right),  \quad R_0>1, \\
    B & = B_1 N + B_2 + \text{O}\left(\frac{1}{N}\right), \quad R_0>1, \\ 
    C & = C_1 N + \text{O}(1), \quad R_0>1,\\
    D & = D_1 N + \text{O}(1), \quad R_0>1,\\
    E & = E_1 N + \text{O}(1), \quad R_0>1,
\end{align}
where
 \begin{align}
   A_1 & = 1 - \alpha R_0 x_1 x_3 - x_1, \\
   A_2 & = - \alpha R_0 (x_1 x_4 + x_2 x_3 + x_6) - x2, \\ 
   B_1 & = \alpha R_0 x_1 x_3 - \alpha x_3, \\ 
   B_2 & = \alpha R_0 (x_1 x_4 + x_2 x_3 + x_6) - \alpha x_4, \\  
   C_1 & = \alpha R_0 x_1 x_3 - 2 x_1 x_6 - 2 x_3 x_5) + 1 + x_1 - 2 x_5,  \\
   D_1 & = \alpha R_0 (x_1 x_6 + x_3 x_5 - x_1 x_3 - x_1 x_7 -x_3 x_6) - (\alpha + 1) x_6, \\
   E_1 & = \alpha R_0 (x_1 x_3 + 2 x_1 x_7 + 2 x_3 x_6) + \alpha (x_3 - 2 x_7). 
\end{align}

The 7 equations formed by setting each of these expressions equal to zero can be solved for the 
7 quantities $x_1$--$x_7$. 
The equation $B_1 = 0$ has two roots, namely $x_1 = 1/R_0$ and $x_3=0$. 
The second of these roots is excluded as the only spurious solution that appears in this method. 
The non-spurious results are as follows:
\begin{align}
    x_1 & = \frac{1}{R_0}, \\
    x_2 & = \frac{\alpha}{R_0-1}, \\ 
    x_3 & = \frac{R_0-1}{\alpha R_0}, \\ 
    x_4 & = - \frac{1}{R_0-1}, \\ 
    x_5 & =  \frac{R_0+\alpha}{R_0^2}, \\ 
    x_6 & = - \frac{1}{R_0}, \\ 
    x_7 & = \frac{R_0^2 + \alpha (R_0-1)}{\alpha R_0^2}.
\end{align}

We conclude that asymptotic approximations of the 5 cumulants of orders 1 and 2 can be written 
as follows:
\begin{align}
   \kappa_{10} & = \frac{1}{R_0}N + \frac{\alpha}{R_0-1} + \text{O}\left(\frac{1}{N}\right), 
       \quad R_0>1, \\
   \kappa_{01} & = \frac{R_0-1}{\alpha R_0}N - \frac{1}{R_0-1} + \text{O}\left(\frac{1}{N}\right), 
       \quad R_0>1, \\
   \kappa_{20} & = \frac{R_0+\alpha}{R_0^2}N + \text{O}(1), \quad R_0>1, \\ 
    \kappa_{11} & = - \frac{1}{R_0}N + \text{O}(1), \quad R_0>1, \\ 
    \kappa_{02} & = \frac{R_0^2+\alpha(R_0-1)}{\alpha R_0^2}N + \text{O}(1), \quad R_0>1. 
\end{align}

We note that the one-term approximations of the cumulants of orders 1 and 2 of the 
bivariate QSD coincide with the approximations derived by the aid of diffusion approximation and 
given in \eqref{5.3}--\eqref{5.7}. 
In this section we have derived two-term approximations of the cumulants of order 1, and 
one-term approximations of the cumulants of order 2. 
The work was based on 5 ODEs for the cumulants of orders 1 and 2, and it required us to determine 
7 coefficients. 
Improved approximations of the cumulants of the first two orders can be achieved by increasing 
the number of ODEs of cumulants that are anahyzed. 
By including cumulants up to order 3, we are led to analyze ODEs of 9 cumulants, and to determine 
a total of 16 coefficients.  
(We note that the approximation of the expected number of infected individuals in QSD, which equals  
$\kappa_{01}$, has an incorrect sign in front of the second term in (14.10) in N{\aa}sell (2005).)

\section{Concluding comments}

The alternative to the cumulant closure method that we have presented here has 
several properties that make it attractive from a modelling standpoint. 
One such property is that the condition for validity of the approximations that it leads to 
($R_0>1$ and $N$ large both for the univariate stochastic logistic Verhulst model and for 
the bivariate stochastic SIR model) can be specified with the aid of the parameters of the model. 
A second attractive property of the new method is that the magnitude of the error that it  
causes can easily be evaluated, as is the case for any asymptotic approximation. 
A third advantage is that the nuisance of dealing with spurious solutions is greatly simplified. 

As a fourth advantage we note that the forms of the approximations produced by the new method are pleasingly simple in comparison with the approximations from the cumulant closure method. 
As seen by \eqref{4.20}--\eqref{4.22}, the approximations from the new method 
of the first 3 cumulants for the Verhulst model are given as explicit functions of the model parameters. 
In similarity to this, Matis and Kiffe (1996) report that closed form approximations are available 
for the first 3 cumulants that they derive using moment closure. 
However, they describe them as too complex to be useful in practical applications.  
Therefore, they do not present them, but they use them for numerical evaluations.

A comparison between the forms of the approximations of the two approaches is, however, possible 
in the simpler case when approximations of the first 2 cumulants are considered,  based on 
a study of the ODEs for $\kappa_1$ and $\kappa_2$. 
The results reported by Matis and Kiffe (1996) for the Verhulst model take the following forms, using the 
second parameter space described in Section 2: 
\begin{align} \label{6.1}
   \kappa_1 & \approx \frac{3(R_0-1)+\gamma_1/\mu}{4(R_0+\alpha)}N, \\ \label{6.2}
   \kappa_2 & \approx 
     \frac{(R_0-1)^2 + 4(\alpha+1)R_0/N - (R_0-1)\gamma_1/\mu}{8(R_0+\alpha)^2} N^2,
\end{align}
where 
\begin{equation} \label{6.3}
   \gamma_1 = \mu [(R_0-1)^2 - 8(\alpha+1)R_0/N]^{1/2}.
\end{equation}
The corresponding approximations of $\kappa_1$ and $\kappa_2$ derived with the new method in
this case are found from \eqref{4.20} and \eqref{4.21} to be as follows: 
\begin{align}
   \kappa_1 & = \frac{R_0-1}{R_0+\alpha}N + \frac{(1+\alpha) R_0}{(R_0+\alpha)(R_0-1} 
      + \text{O}\left(\frac{1}{N}\right), \quad R_0>1, \\
   \kappa_2 & = \frac{(1+\alpha)R_0}{(R_0+\alpha)^2}N + \text{O}(1), \quad R_0>1.
\end{align}
Both sets of approximations are given as functions of the parameters, but the new results are 
simpler and more attractive in form and easier to work with.

The number of terms in the asymptotic approximations of the first 3 stationary cumulants 
of the QSD for the stochastic logistic Verhulst model is seen by \eqref{4.20}--\eqref{4.22} to vary. 
It depends on the order of the cumulant, and also on the number of ODEs that are analyzed. 
Thus, if we start with 1 ODE (for $\kappa_1$), then we can use it to derive 1 term in the 
asymptotic approximation of $\kappa_1$. 
By using 2 ODEs (for $\kappa_1$ and $\kappa_2$), we can derive 2 terms in the 
asymptotic approximation of $\kappa_1$ and 1 term in the asymptotic approximation 
of $\kappa_2$. 
Next, by using 3 ODEs (for $\kappa_1$, $\kappa_2$, and $\kappa_3$), 
as was done in Section 4, we are led to 3, 2, and 1 terms, respectively, in the asymptotic
approximations of $\kappa_1$, $\kappa_2$, and $\kappa_3$. 
It is easy to conjecture extensions to analyses based on more than 3 ODEs, and also extensions 
to bivariate and multivariate models. 
So by increasing the number of ODEs whose critical points are determined, we increase the 
information about the asymptotic behaviors of the low order cumulants. 
It is possible to derive asymptotic approximations of cumulants of higher order. 
Such an extension would at the same time give additional terms of the 
asymptotic approximations of the cumulants of lower order.

Our results place three requirements on the modelling work that are different 
from early applications of moment closure. 
The first one is that  the model parametrization should be similar to the second one described in 
Section 2, where the maximum population size $N$ and the basic reproduction 
ratio $R_0$ appear as explicit parameters. 
The second requirement deals with the situation when the unconditioned random variable 
or random vector has an absorbing state or an absorbing class of states. 
The analysis should then deal with the random variable or vector conditioned on 
non-absorption. 
The third requirement is that orders of magnitude of the first few cumulants need to be determined
before the forms of the asymptotic approximations can be formulated.

We note that the procedure that we describe here as an alternative to moment 
closure gives useful results only in the parameter region where $R_0>1$. 
The ODEs \eqref{3.7}--\eqref{3.9} for the stochastic logistic Verhulst model are clearly valid in the remaining two parameter regions, namely where $R_0<1$, and where $\rho$ defined by 
\eqref{2.7} is constant. 
However, the step required to formulate the right-hand sides of the ODEs for the first three 
cumulants in these two regions requires more information about $q_1$ than what was 
needed in the parameter region $R_0>1$, since $q_1$ is not exponentially small in these regions. 
Additonal studies are required to supply such information.

The main steps that need to be taken to apply the method that we have described here 
to other models that have been analyzed with the aid of a moment closure method 
should be clear from our description: 
The ODEs for the first few cumulants are the same as those that have been derived
in the moment closure method. 
A reparameterization may be necessary so that one can identify one parameter that takes large values. 
In addition one needs to determine the orders of magnitude of the cumulants that appear in 
these ODEs. 
We have shown that both numerical evaluations and diffusion approximations can be used 
to give a guidance in this regard.

The advantages of the method for deriving asymptotic approximations 
described here are not limited to the two specific stochastic models that we have studied. 
We expect that the method is equally useful on any continuous-time discrete-state Markov 
Chain with nonlinear transition rates. 
This holds in particular for the open problems formed by the stochastic 
versions of some of the deterministic population models studied by Tsoularis 
and Wallace (2002). 
Similar considerations may be valid in some of the many mathematical models treated in the 
review of moment closure methods given by Kuehn (2016).

Kuehn (2016) states that the system of ODEs must be closed in order to make it tractable 
for analytic or numerical techniques. 
This tenet has served as a basic principle in all work where moment closure has played a role. 
However, our approach shows that such closures are not necessary if one is satisfied with approximate results. 
We note furthermore that the assumptions that are introduced in order to achieve moment closure invariably lead to approximations. 
This means that the requirement of exactness of solutions must be abandoned. 
A consequence of this is that we can also abandon the requirement of finding a closure of the 
system of ODEs, and in this way avail ourselves of the power of asymptotic approximations.

In his review of moment closure methods, Kuehn (2016) remarks that moment closure approximations 
work well in practice, but that it has been difficult to justify them rigorously.  
He uses geometric invariant manifold theory to argue that it may be difficult to determine 
the approximation errors.

We have shown in this paper that fundamental mistakes were made in early studies of cumulants of the 
stochastic logistic Verhulst model. 
It is remarkable that these mistakes have remained undiscovered until 57 years have passed 
since the publication of the 1960 paper by Bartlett, Gower, and Leslie.

It appears that most authors dealing with moment closure have been content with the results, 
since numerical evaluations indicate that they are reasonable, when one excludes those cases 
where the approximations do not work at all. 
Even so, it is likely that many adherents of moment closure have been aware of some of 
the methodological weaknesses with these methods, and that therefore it has been natural for 
them to search for improvements.  
In view of this, it is noteworthy that it took 60 years after Whittle's basic paper (1957) was published before the improved methods described here were discovered.

\end{document}